\documentclass[journal]{IEEEtran}

\usepackage{cite}      
\usepackage{subfigure} 
\usepackage{url}       
\usepackage{algorithm}
\usepackage{array}
\usepackage{amsmath,amsfonts,amssymb,graphicx}
\usepackage{bm}
\usepackage{mathrsfs}
\usepackage{dsfont}
\usepackage{multicol}
\usepackage[end]{algpseudocode}
\usepackage{multirow}

\newcommand{\commnt}[1] {$//$ \textsc{#1} }
\newtheorem{theorem}{\textbf{Theorem}}
\newtheorem{lemma}{\textbf{Lemma}}

\newtheorem{fact}{\textbf{Fact}}
\newtheorem{proof}{\textbf{Proof}}

\begin{document}
\title{The Non-Bayesian Restless Multi-Armed Bandit: A Case of Near-Logarithmic Strict Regret\thanks{This research was sponsored in part by the
U.S. Army Research Laboratory under the Network Science
Collaborative Technology Alliance, Agreement Number
W911NF-09-2-0053. The work of Q. Zhao was supported by the Army Research Office under Grant W911NF-08-1-0467. This is
an extended, full version of a paper that appeared in ICASSP
2011~\cite{ICASSPpaper}.}}

\author{\IEEEauthorblockN{Wenhan Dai$^\dag$, Yi Gai$^\ddag$, Bhaskar
Krishnamachari$^\ddag$ and Qing Zhao$^\S$}
\IEEEauthorblockA{\\$^\dag$ Department of Aeronautics and Astronautics\\
Massachusetts Institute of Technology, Cambridge, MA 02139, USA\\
$^\ddag$ Ming Hsieh Department of Electrical Engineering\\
University of Southern California, Los Angeles, CA 90089, USA\\
$^\S$Department of Electrical and Computer Engineering\\
University of California, Davis, CA 95616, USA\\
 Email: whdai@mit.edu; $\{$ygai,bkrishna$\}$@usc.edu; qzhao@ucdavis.edu
 }
 }

\markboth{}{}
\maketitle

\begin{abstract}
In the classic Bayesian restless multi-armed bandit (RMAB) problem,
there are $N$ arms, with rewards on all arms evolving at each time
as Markov chains with known parameters. A player seeks to activate
$K \geq 1$ arms at each time in order to maximize the expected total
reward obtained over multiple plays. RMAB is a challenging problem
that is known to be PSPACE-hard in general. We consider in this work
the even harder non-Bayesian RMAB, in which the parameters of the
Markov chain are assumed to be unknown \emph{a priori}. We develop
an original approach to this problem that is applicable when the
corresponding Bayesian problem has the structure that, depending on
the known parameter values, the optimal solution is one of a
prescribed finite set of policies. In such settings, we propose to
learn the optimal policy for the non-Bayesian RMAB by employing a
suitable meta-policy which treats each policy from this finite set
as an arm in a different non-Bayesian multi-armed bandit problem for
which a single-arm selection policy is optimal. We demonstrate this
approach by developing a novel sensing policy for opportunistic
spectrum access over unknown dynamic channels. We prove that our
policy achieves near-logarithmic regret (the difference in expected
reward compared to a model-aware genie), which leads to the same
average reward that can be achieved by the optimal policy under a
known model. This is the first such result in the literature for a
non-Bayesian RMAB. For our proof, we also develop a novel
generalization of the Chernoff-Hoeffding bound.
\end{abstract}


\section{Introduction}
\IEEEPARstart{M}{ulti-armed} bandit (MAB) problems are fundamental
tools for optimal decision making in dynamic, uncertain
environments. In a multi-armed bandit problem, there are $N$ independent arms
each generating stochastic rewards, and a player seeks a policy to
activate $K \geq 1$ arms at each time in order to maximize the
expected total reward obtained over multiple plays. A particularly
challenging variant of these problems is the restless multi-armed
bandit problem (RMAB)~\cite{Whittle}, in which the rewards on all
arms (whether or not they are played) evolve at each time as Markov
chains.

MAB problems can be broadly classified as Bayesian or non-Bayesian.
In a Bayesian MAB, there is a prior distribution on the arm rewards
that is updated based on observations at each step and a
known-parameter model for the evolution of the rewards. In a
non-Bayesian MAB, a probabilistic belief update is not possible
because there is no prior distribution and/or the parameters of the
underlying probabilistic model are unknown. In the case of
non-Bayesian MAB problems, the objective is to design an arm
selection policy that minimizes regret, defined as the gap between
the expected reward that can be achieved by a genie that knows the
parameters, and that obtained by the given policy. It is desirable
to have a regret that grows as slowly as possible over time. In
particular, if the regret is sub-linear, the average regret per slot
tends to zero over time, and the policy achieves the maximum average
reward that can be achieved under a known model.

Even in the Bayesian case, where the parameters of the Markov chains
are known, the restless multi-armed bandit problem is difficult to
solve, and has been proved to be PSPACE hard in general~\cite{Papadimitriou}.
One approach to this problem has been Whittle's index, which is
asymptotically optimal under certain regimes~\cite{WeberWeiss}; however it does not
always exist, and even when it does, it is not easy to compute. It
is only in very recent work that non-trivial tractable classes of
RMAB where Whittle's index exists and is computable have been
identified~\cite{LiuZhaoIndex}~\cite{LeNyDahleh}.

We consider in this work the even harder non-Bayesian RMAB, in which
the parameters of the Markov chain are assumed to be unknown
\emph{a priori}. Our main contribution in this work is a novel
approach to this problem that is applicable when the
corresponding Bayesian RMAB problem has the structure that the
parameter space can be partitioned into a finite number of sets, for
each of which there is a single optimal policy. For RMABs satisfying this finite partition property, our approach
is to develop a meta-policy that treats these policies as
arms in a different non-Bayesian multi-armed bandit problem for
which a single arm selection policy is optimal for the genie, and learn from reward observations which policy from this finite set gives the best
performance.

We demonstrate our approach on a practical problem pertaining to
dynamic spectrum sensing. In this problem, we consider a scenario
where a secondary user must select one of $N$ channels to sense at
each time to maximize its expected reward from transmission
opportunities. If the primary user occupancy on each channel is
modeled as an identical but independent Markov chain with unknown
parameters, we obtain a non-Bayesian RMAB with the requisite
structure. We develop an efficient new multi-channel cognitive
sensing policy for unknown dynamic channels based on the above
approach. We prove for $N=2,3$ that this policy achieves regret (the
gap between the expected optimal reward obtained by a model-aware
genie and that obtained by the given policy) that is bounded
uniformly over time $n$ by a function that grows as $O(G(n)\cdot
\log n)$, where $G(n)$ can be any arbitrarily slowly diverging
non-decreasing sequence. For the general case, this policy achieves
the average reward of the myopic policy which is conjectured, based
on extensive numerical studies, to be optimal for the corresponding
problem with known parameters. This is the first non-Bayesian RMAB
policy that achieves the maximum average reward defined by the
optimal policy under a known model.

\section{Related Work}
\label{sec:Rel}
\subsection{Bayesian MAB}

The Bayesian MAB takes a probabilistic viewpoint toward the system unknown parameters. By treating the player's \emph{a posterior} probabilistic knowledge (updated
from the \emph{a priori} distribution using past observations) on the
unknown parameters as the system state, Bellman in 1956 abstracted and generalized the Bayesain MAB to a special class of Markov decision processes (MDP)~\cite{Bellman:1956}. Specifically, there are N independent arms with fully observable states. One arm is activated at each time, and only the
activated arm changes state as per a known Markov process and offers
a state-dependent reward. This general MDP formulation of the problem naturally leads to a stochastic dynamic programming solution based on backward induction. However, such an approach incurs exponential complexity with respect to the number of arms. The problem of finding a simpler approach remained open till 1972 when Gittins and Jones~\cite{Gittins:1972} presented
a forward-induction approach in which an index is calculated for
each arm depending only on the process of that arm, and the arm with
the highest index at its current state is selected at each time. This result shows that arms can be decoupled when seeking the optimal activation rule, consequently reducing the complexity from exponential to linear in terms of the number of arms. Several researchers have since developed alternative proofs of the optimality of this approach, which has come to be known as the
Gittins-index~\cite{Whittle:1980, Varaiya:1985, Weber:1992,
Berry:1985, Tsitsiklis:1986, Mandelbaum:1986, Ishikida:1986,
Kaspi:1998}. Several variants of the basic classical Bayesian MAB
have been proposed and investigated, including arm-acquiring
bandits~\cite{Whittle:1981}, superprocess bandits~\cite{Nash:1973,
Whittle:1980}, bandits with switching penalties~\cite{Banks:1994,
Asawa:1996}, and multiple simultaneous plays~\cite{Agrawal:1990, Pandelis:1999}.

A particularly important variant of the classic MAB is the
\emph{restless bandit problem} posed by Whittle in
1988~\cite{Whittle:1988}, in which the passive arms also change
state (to model system dynamics that cannot be directly controlled).
The structure of the optimal solution for this problem in general
remains unknown, and has been shown to be PSPACE-hard by
Papadimitriou and Tsitsiklis~\cite{Papadimitriou:1999}. Whittle
proposed an index policy for this problem that is optimal under a
relaxed constraint of an average number of arms played as well as
asymptotically under certain conditions~\cite{Weber:1990}; for
many problems, this Whittle-index policy has numerically been found
to offer near-optimal performance. However, Whittle index is not
guaranteed to exist. Its existence (the so-called indexability) is
difficult to check, and the index can be computationally expensive to calculate when it does exist. General analytical results on the optimality of Whittle index in the finite regime have also eluded the research
community up to today. There are numerical approaches for testing indexability and calculating Whittle index (see, for example, \cite{Nino,Nino:07TOP}). Constant-factor approximation algorithms for restless bandits have also been explored in the literature~\cite{Guha,Guha1}.

Among the recent work that contributes to the fundamental
understanding of the basic structure of the optimal policies for a
class of restless bandits with known models, myopic
policy~\cite{zhao:krishnamachari:twc2008,
javidi:krishnamachari:icc2008, ahmad:liu:it2009} has a simple
semi-universal round-robin structure. It has been shown  that the
myopic policy is optimal for $N=2,3$, and for any $N$ in the case of
positively correlated channels. The optimality of the myopic policy
for $N>3$ negatively correlated channels is conjectured for the
infinite-horizon case. Our work provides the first efficient
solution to the non-Bayesian version of this class of problems,
making use of the semi-universal structure identified in
\cite{zhao:krishnamachari:twc2008}.

\subsection{Non-Bayesian MAB}

In the non-Bayesian formulation of MAB, the unknown parameters in the system dynamics are treated as deterministic quantities; no \emph{a priori} probabilistic knowledge about the unknowns is required. The basic form
of the problem is the optimal sequential activation of $N$ independent arms, each associated with an i.i.d. reward process with an unknown mean. The performance of an arm activation policy is measured by \emph{regret} (also known as the cost of learning) defined as the difference between
the total expected reward that could be obtained by an omniscient
player that knows the parameters of the reward model and the policy
in question (which has to learn these parameters through statistical
observations). Notice that with the reward model known, the omniscient player will always activate the arm with the highest reward mean. The essence of the problem is thus to identify the best arm without exploring the bad arms too often in order to minimize the regret. In particular, it is
desirable to have a sub-linear regret function with respect to time,
as under this condition the time-averaged regret goes to zero, and the
slower the regret growth rate, the faster the system converges to the
same maximum average reward achievable under the known-model case.

Lai and Robbins~\cite{Lai:Robbins} proved in 1985 that the lower
bound of regret is logarithmic in time, and proposed the first
policy that achieved the optimal logarithmic regret for non-Bayesian
MABs in which the rewards are i.i.d over time and obtained from a
distribution that can be characterized by a single-parameter.
Anantharam \emph{et al.} extended this result to multiple
simultaneous arm plays, as well as single-parameter Markovian
\emph{rested}
rewards~\cite{Anantharam:1987:part1,Anantharam:1987:part2}. Other
policies achieving logarithmic regret under different assumptions
about the i.i.d. reward model have been developed by
Agrawal~\cite{Agrawal:1995} and Auer~\emph{et al.}~\cite{Auer:2002}.
In particular, Auer~\emph{et al.}'s UCB1 policy applies to i.i.d.
reward distributions with finite support, and achieves logarithmic
regret with a known leading constant uniformly bounded over time.

The focus of this paper is on the non-Bayesian RMAB. There are two
parallel investigations on non-Bayesian RMAB problems given in
\cite{TekinLiu,LiuLiuZhao}, where a more general RMAB model is
considered but under a much weaker definition of regret.
Specifically, in \cite{TekinLiu,LiuLiuZhao}, regret is defined with
respect to the maximum reward that can be offered by a \emph{single
arm/channel}. Note that for RMAB with a known model, staying with
the best arm is suboptimal in general. Thus, a sublinear regret
under this definition does not imply the maximum average reward, and
the deviation from the maximum average reward can be arbitrarily
large. In contrast to these works, this paper shows sublinear regret
with respect to the maximum reward that can be obtained by the
optimal policy played by a genie that knows the underlying
transition matrix.

\section{A New Approach for non-Bayesian RMAB}
In multi-arm bandit problem, there are multiple arms and each of
them yields a stochastic reward when played. The player sequentially picks
one arm at each time, aiming to maximize the total expected reward collected over
time. If the rewards on all arms
are modeled as Markov chains and all arms always keep activated
whether they are selected, it is classified as restless multi-armed
bandit problem (RMAB). In Bayesian RMAB, the parameters of the
Markov chain are known and in non-Baysian RMAB, the model for the
reward process is a priori unknown to the user.

We first describe a structured class of finite-option Bayesian RMAB
problems that we will refer to as $\Psi_m$. Let $\mathcal{B}(P)$ be
a Bayesian RMAB problem with the Markovian evolution of arms
described by the transition matrix $P$. We say that $\mathcal{B}(P)
\in \Psi_m$ if and only if there exists a partition of the parameter
values $P$ into a finite number of $m$ sets $\{ S_1, S_2, ... S_m\}
$ and a set of policies $\pi_i$ $(\forall i = 1\ldots m)$ with $\pi_i$ being optimal whenever $P \in S_i$.
Despite the general hardness of the RMAB problem, problems with such
structure do exist, as has been shown
in~\cite{zhao:krishnamachari:twc2008, ahmad:liu:it2009,
LiuZhaoIndex}.

We propose a solution to the non-Bayesian version of the problem
that leverages the finite solution option structure of the corresponding Bayesian version ($\mathcal{B}(P) \in \Psi_m$).
In this case, although the player does not know the exact parameter
$P$, it must be true that one of the $m$ policies $\pi_i$ will yield
the highest expected reward (corresponding to the set $S_i$ that
contains the true, unknown $P$). These policies can thus be treated
as arms in a different non-Bayesian multi-armed bandit problem for
which a single-arm selection policy is optimal for the genie. Then,
a suitable meta-policy that sequentially operates these policies
while trying to minimize regret can be adopted. This can be done
with an algorithm based on the well-known schemes proposed by Lai
and Robbins~\cite{Lai:Robbins}, and Auer \emph{et
al}~\cite{Auer:2002}.

One subtle issue that must be handled in adopting such an algorithm
as a meta-policy is how long to play each policy. An ideal constant
length of play could be determined only with knowledge of the
underlying unknown parameters $P$. To circumvent this difficulty, our approach is to have the
duration for which each policy is operated slowly increase over
time.

In the following, we demonstrate this novel meta-policy approach
using the dynamic spectrum access problem discussed
in~\cite{zhao:krishnamachari:twc2008, ahmad:liu:it2009} where the
Bayesian version of the RMAB has been shown to belong to the class
$\Psi_2$. For this problem, we show that our approach yields a
policy with provably near-logarithmic regret, thus achieving the
same average reward offered by the optimal RMAB policy under a known
model.

\section{Dynamic Spectrum Access under Unknown Models}
\label{sec:problem}

We consider a slotted system where a secondary user is trying to
access $N$ independent channels, with the availability of each
channel evolving as a two-state Markov chain with identical
transition matrix $\mathbf{P}$ that is \emph{a priori} unknown to
the user. The user can only see the state of the sensed channel. If
the user selects channel $i$ at time $t$, and upon sensing finds the
state of the channel $S_i(t)$ to be 1, it receives a unit reward for
transmitting. If it instead finds the channel to be busy, i.e.,
$S_i(t) = 0$, it gets no reward at that time. The user aims to
maximize its expected total reward (throughput) over some time
horizon by choosing judiciously a sensing policy that governs the
channel selection in each slot. We are interested in designing
policies that perform well with respect to \emph{regret}, which is
defined as the difference between the expected reward that could be
obtained using the omniscient policy $\pi^*$ that knows the
transition matrix $\mathbf{P}$, and that obtained by the given
policy $\pi$. The regret at time $n$ can be expressed as:
\begin{equation}
\begin{split}
 R(\mathbf{P},\Omega(1),n) & =
 \mathbb{E}^{\pi^*}[\Sigma_{t=1}^{n}Y^{\pi^*}(\mathbf{P},\Omega(1),t)]
 \\
 & -
 \mathbb{E}^{\pi}[\Sigma_{t=1}^{n}Y^{\pi}(\mathbf{P},\Omega(1),t)],
\end{split}
\end{equation}
where $\omega_i(1)$ is the initial probability that $S_i(1)=1$, $\mathbf{P}$ is the transition matrix
of each channel, $Y^{\pi^*}(\mathbf{P},\Omega(1),t)] $ is the reward obtained in time $t$ with the
optimal policy, $ Y^{\pi}(\mathbf{P},\Omega(1),t)$ is the reward obtained in time $t$ with the
given policy. We denote $\Omega(t)\triangleq [\omega_1(t), \ldots, \omega_N(t)]$ as the belief
vector where $\omega_i(t)$ is the conditional probability that $S_i(t) = 1$ (and let
$\Omega(1)=[\omega_1(1),\ldots,\omega_N(1)]$ denote the initial belief vector used in the myopic
sensing algorithm~\cite{zhao:krishnamachari:twc2008}).

\section{Policy Construction}
\label{sec:algorithm}

As has been shown in~\cite{zhao:krishnamachari:twc2008}, the myopic policy has a
simple structure for switching between channels that depends only on
the correlation sign of the transition matrix $\mathbf{P}$, i.e.
whether $p_{11}\geq p_{01}$ (positively correlated) or $p_{11}<
p_{01}$ (negatively correlated).

In particular, if the channel is positively correlated, then the
myopic policy corresponds to
\begin{itemize}
\item \textbf{Policy $\pi_1$}: stay on a channel whenever it shows a
``1''
and switch on a ``0'' to the channel visited the longest ago.
\end{itemize}

If the channel is negatively correlated, then it corresponds to
\begin{itemize}
\item \textbf{Policy $\pi_2$:} staying on a channel when it shows a ``0'', and
switching as soon as ``1'' is observed, to either the channel most
recently visited among those visited an even number of steps before,
or if there are no such channels, to the one visited the longest
ago.
\end{itemize}

To be more specific, we give the structure of myopic sensing~\cite{zhao:krishnamachari:twc2008}. In myopic sensing, the concept \emph{circular order} is very important. A circular order $\kappa = (n_1,n_2,\cdots,n_N)$ is equivalent to $(n_i,n_{i+1},\cdots,n_N,n_1,n_2,\cdots,n_{i-1})$ for any $1 \leq i \leq N$. For a circular order $\kappa$, denote $-\kappa$ as its reverse circular order. For a channel $i$, denote $i_\kappa^+$ as its next channel in the circular order $\kappa$. With these notations, we present the structure of the myopic sensing.

Let $\Omega(1) = [\omega_1(1),\cdots,\omega_N(1)]$ denote the initial belief vector. The circular order $\kappa(1)$ in time slot 1 depends on the order of $\Omega(1)$: $\kappa(1) = (n_1,n_2,\cdots,n_N)$ implies that $\omega_{n_1}(1) \leq \omega_{n_2}(1) \leq \cdots \leq \omega_{n_N}(1)$. Let $\hat{a}(t)$ denote myopic action in time $t$. We have $\hat{a}(1) = \arg\max_{i=1,2,\cdots,N}\omega_i(1)$ and for $t > 1$, the myopic action $\hat{a}(t)$ is given as follows.

\begin{itemize}
\item \textbf{Policy $\pi_1$}($p_{11}>p_{01}$):
\begin{equation} \nonumber
\hat{a}(t)=\left\{ \begin{aligned}
         \hat{a}(t-1),   if \mathbf{\emph{S}}_{\hat{a}(t-1)}(t-1)=1 \\
                  \hat{a}(t-1)_{\kappa(t)}^+,   if S_{\hat{a}(t-1)}(t-1)=0
                          \end{aligned} \right.
\end{equation}
where $\kappa(t) \equiv \kappa(1)$.
\end{itemize}

\begin{itemize}
\item \textbf{Policy $\pi_2$}($p_{11}<p_{01}$):
\begin{equation} \nonumber
\hat{a}(t)=\left\{ \begin{aligned}
         \hat{a}(t-1),   if \mathbf{\emph{S}}_{\hat{a}(t-1)}(t-1)=0 \\
                  \hat{a}(t-1)_{\kappa(t)}^+,   if S_{\hat{a}(t-1)}(t-1)=1
                          \end{aligned} \right.
\end{equation}
where $\kappa(t) = \kappa(1)$ when $t$ is odd and $\kappa(t) = -\kappa(1)$ when $t$ is even.
\end{itemize}

Furthermore, as mentioned in section \ref{sec:Rel}, it has been
shown in~\cite{zhao:krishnamachari:twc2008, ahmad:liu:it2009} that
the myopic policy is optimal for $N=2,3$. As a consequence, this
special class of RMAB has the required finite dependence on its
model as described in section \ref{sec:problem}; specifically, it
belongs to $\Psi_2$. We can thus apply the general approach
described in section \ref{sec:problem}. Specifically, the algorithm
treats these two policies as arms in a classic non-Bayesian
multi-armed bandit problem, with the goal of learning which one
gives the higher reward.

A key question is how long to operate each arm at each step. The
analysis we present in the next section shows that it is desirable
to slowly increase the duration of each step using any (arbitrarily
slowly) divergent non-decreasing sequence of positive integers
$\{K_n\}_{n=1}^\infty$.

The channel sensing policy we thus construct is shown in Algorithm
\ref{alg:sensing}, in which we use the UCB1 policy proposed by Auer
\emph{et al.} in ~\cite{Auer:2002} as the meta-policy.

\begin{algorithm} [ht]
\caption{Sensing Policy for Unknown Dynamic Channels}
\label{alg:sensing}
\begin{algorithmic}[1]
\State \commnt{ Initialization}

\State Let $\{K_n\}_{n=1}^\infty$ be any arbitrarily slowly
divergent non-decreasing sequence of positive integers.

\State Play policy $\pi_1$ for $K_1$ times, denote $\hat{A}_1$ as
the sample mean of these $K_1$ rewards \State Play policy $\pi_2$
for $K_2$ times, denote $\hat{A}_2$ as the sample mean of these
$K_2$ rewards \State $\hat{X}_1 = \hat{A}_1$, $\hat{X}_2 =
\hat{A}_2$ \State $n = K_1+K_2$ \State $i=3$, $i_1=1$, $i_2=1$

\State \commnt{Main loop}

\While {1}
    \State Find $j \in \{ 1, 2\}$ such that
    $j = \arg\max\frac{\hat{X}_j}{i_j}+\sqrt{\frac{L\ln{n}}{i_j}}$ \label{line:9}(L can be any constant larger than 2)
    \State $i_j = i_j+1$
    \State Play policy $\pi_j$ for $K_i$ times, let $\hat{A}_{j}(i_j)$ record the sample mean of these $K_i$ rewards
    \State $\hat{X}_j = \hat{X}_j + \hat{A}_{j}(i_j)$
    \State $i = i+1$
    \State $n = n + K_i$;
\EndWhile
\end{algorithmic}
\end{algorithm}

\section{Regret Analysis}
\label{sec:regret}

We first define the discrete function $G(n)$ which represents the
value of $K_i$ at the $n^{th}$ time step in
Algorithm~\ref{alg:sensing}:
\begin{equation}G(n) = \min\limits_I K_I~~s.t.~~\sum\limits_{i=1}^I K_i \geq  n
\end{equation}
Note that since $K_i$ can be any arbitrarily slow non-decreasing
diverging sequence, $G(n)$ can also grow arbitrarily slowly.

The following theorem states that the regret of our policy grows
close to logarithmically with time.

\begin{theorem} \label{theorem:regret} For the dynamic spectrum access
problem with $N=2,3$ i.i.d. channels with unknown transition matrix
$\mathbf{P}$, the expected regret with Algorithm \ref{alg:sensing}
after $n$ time steps is at most $Z_1 G(n)\ln(n) + Z_2 \ln(n) + Z_3
G(n) + Z_4$, where $Z_1,Z_2,Z_3,Z_4$ are constants only related to
$\mathbf{P}$.
\end{theorem}

The proof of Theorem \ref{theorem:regret} uses one fact and two lemmas.

\fact \label{fact:chernoff} (Chernoff-Hoeffding
bound~\cite{Pollard}) Let $X_1,\cdots,X_n$ be random variables with
common range $[0,1]$ such that
$\mathbb{E}[X_t|X_1,\cdots,X_{t-1}]=\mu$. Let $S_n =
X_1+\cdots+X_n$. Then for all a $\geq 0$
\begin{equation}
\mathbb{P}\{S_n \geq n\mu + a\}\leq e^{-2a^2/n}; \mathbb{P}\{S_n \leq
n\mu - a\}\leq e^{-2a^2/n}
\end{equation}

Our first lemma is a non-trivial generalization of the
Chernoff-Hoeffding bound, that allows for bounded differences
between the conditional expectations of a sequence of random variables
that are revealed sequentially:

\lemma \label{lemma:chernoff} Let $X_1,\cdots,X_n$ be random
variables with range $[0,b]$ and such that
$|\mathbb{E}[X_t|X_1,\cdots,X_{t-1}]- \mu | \leq C$. $C$ is a
constant number such that $0< C < \mu$. Let $S_n =
X_1+\cdots+X_n$. Then for all  $a \geq 0$,
\begin{equation} \label{eqn:chernoff1}
\mathbb{P}\{S_n \geq n(\mu+C) + a\}\leq
e^{-2(\frac{a(\mu-C)}{b(\mu+C)})^2/n}
\end{equation} and
\begin{equation} \label{eqn:chernoff2}
\mathbb{P}\{S_n \leq n(\mu-C) - a\}\leq e^{-2(a/b)^2/n}
\end{equation}

\begin{proof}

We first prove (\ref{eqn:chernoff1}). We generate random variables
$\hat{X}_1, \hat{X}_2,\cdots, \hat{X}_n$ as follows:

$\hat{X}_1 = (\mu+C)\frac{X_1}{\mathbb{E}[X_1]},$

$\hat{X}_2 = (\mu+C)\frac{X_2}{\mathbb{E}[X_2|\hat{X}_1]},$

$\cdots$

$\hat{X}_t =
(\mu+C)\frac{X_t}{\mathbb{E}[X_t|\hat{X}_1,\hat{X}_2,\cdots,\hat{X}_{t-1}]}.$

Note that $|\mathbb{E}[X_t|X_1,\cdots,X_{t-1}]- \mu | \leq C$. Therefore
$|\mathbb{E}[X_t|\hat{X}_1,\cdots,\hat{X}_{t-1}]- \mu | \leq C$
also stands. Hence $\frac{\hat{X}_t}{X_t}$ is at least 1, at most
$\frac{\mu+C}{\mu-C}$. Therefore $\hat{X}_1, \hat{X}_2,\cdots,
\hat{X}_n$ have finite support (they are in the range $[0,
b\frac{\mu+C}{\mu-C}]$). Besides,
$\mathbb{E}[\hat{X}_t|\hat{X}_1,\cdots,\hat{X}_{t-1}] = \mu + C$,
$\forall t$.

Let $\hat{S}_n = \hat{X}_1+\hat{X}_2+\cdots+\hat{X}_n$, then for all
$a \geq 0$,
\begin{equation}
\begin{split}
\mathbb{P}\{S_n \geq n(\mu+C) + a\} & \leq \mathbb{P}\{\hat{S}_n \geq
n(\mu+C) + a\} \\
& \leq e^{-2(\frac{a(\mu-C)}{b(\mu+C)})^2/n}.
\end{split}
\end{equation}
The first inequality stands because $\frac{\hat{X}_t}{X_t} \geq 1$,$\forall t$.
The second inequality stands because of Fact 1.

The proof of (\ref{eqn:chernoff2}) is similar. We generate random
variables $\hat{X}_1',\hat{X}_2',\cdots,\hat{X}_n'$ as follows:

$\hat{X}_1' = (\mu-C)\frac{X_1}{\mathbb{E}[X_1]},$

$\cdots$

$\hat{X}_n' =
(\mu-C)\frac{X_n}{\mathbb{E}[X_n|\hat{X}_1',\hat{X}_2',\cdots,\hat{X}_{n-1}']}.$

Note that $|\mathbb{E}[X_t|X_1,\cdots,X_{t-1}]-\mu|\leq C$,
$|\mathbb{E}[X_t'|\hat{X}_1',\cdots,\hat{X}_{t-1}']- \mu | \leq
C$ also stands. So $\frac{\hat{X}_t'}{X_t}$ is at most 1, at least
$\frac{\mu-C}{\mu+C}$. Therefore, $\hat{X}_1, \hat{X}_2,\cdots,
\hat{X}_n$ have finite support(They are in the range $[0,b]$).
Besides, $\mathbb{E}[\hat{X}_t'|\hat{X}_1',\cdots,\hat{X}_{t-1}'] =
\mu - C$, $\forall t$.

Let $\hat{S}_n' = \hat{X}_1'+\hat{X}_2'+\cdots+\hat{X}_n'$, then for
all $a \geq 0$,
\begin{equation}
\begin{split}
\mathbb{P}\{S_n \leq n(\mu-C) - a\} & \leq \mathbb{P}\{\hat{S}_n' \leq
n(\mu-C) - a\}\\
& \leq e^{-2(a/b)^2/n}.
\end{split}
\end{equation}
The first inequality stands because $\frac{\hat{X}_t'}{X_t} \leq 1$,$\forall t$.
The second inequality stands because of Fact 1.
\end{proof}




The second lemma states that the expected loss of reward for either
policy due to starting with an arbitrary initial belief vector
compared to the reward $U_i(\mathbf{P})$ that would obtained by playing the
policy at steady state is bounded by a constant $C_i(\mathbf{P})$ that
depends only on the policy used and the transition matrix.

Specifically, denote
\begin{equation}
U_i(\mathbf{P},\Omega(1))\triangleq \lim_{T \to
\infty}\frac{\mathbb{E}^{\pi^*}[\mathop{\Sigma}_{t=1}^{T} Y^{\pi^*}(\mathbf{P},\Omega(1),t)]}{T}.
\end{equation}
From the previous work in ~\cite{LiuZhaoIndex}, we know that the limit
above exists and the steady-state throughput $U_i(\mathbf{P},\Omega(1))$ is independent of
the initial belief value $\Omega(1)$. So, we can
rewrite $U_i(\mathbf{P},\Omega(1))$ as $U_i(\mathbf{P})$, to denote the average expected
reward with transition matrix $\mathbf{P}$ using policy $\pi_i (i=1,2)$.

\lemma \label{lemma:2} For any initial belief vector $\Omega(1)$
and any positive integer $M$, if we use policy $\pi_i$
($i=1,2$) for $M$ times, and the summed expectation of the rewards
for these $M$ steps is denoted as
$\mathbb{E}^{\pi_i}[\Sigma_{t=1}^{M}Y^{\pi_i}(\mathbf{P},\Omega(1),t)]$, then
\begin{equation}
|\mathbb{E}^{\pi_i}[\Sigma_{t=1}^{M}Y^{\pi_i}(\mathbf{P},\Omega(1),t)]-M
\cdot U_i(P)|<C_i(P)
\end{equation}

\begin{proof}

As described in section \ref{sec:problem}, let $\kappa(t) = (n_1,
n_2, \cdots, n_N)$ $(n_i\in\{1,2,\dots,N\}, \forall i)$ be the
circular order of channels in slot t, we then have an ordered
channel states $\vec{S(t)} \triangleq [S_{n_1}(t),S_{n_2}(t),\cdots,
S_{n_N}(t)]$. In fact, after sensing at time $t-1$ , channel
sequence $\{n_1, n_2, \cdots, n_N\}$ has a non-decreasing
probability of being in state 1, i.e. $\omega_{n_1}(t) \geq
\omega_{n_2}(t) \geq \cdots \omega_{n_N}(t)$.

The status $\vec{S(t)}$ form a $2^N$-state Markov chain with transition probabilities $\{Q_{\vec{i},\vec{j}}\}$ shown as follows:

for policy $\pi_1$,
\begin{equation}
Q_{\vec{i},\vec{j}}=\left\{ \begin{aligned}
         \prod_{k=1}^N p_{i_k,j_k},   \quad \text{if } i_1 = 1 \\
                  p_{i_1,j_N}\prod_{k=2}^Np_{i_k,j_{k-1}},   \quad \text{if } i_1 = 0
                          \end{aligned} \right.
\end{equation}

and for policy $\pi_2$,
\begin{equation}
Q_{\vec{i},\vec{j}}=\left\{ \begin{aligned}
         \prod_{k=1}^N p_{i_k,j_{N-k+1}},   \quad \text{if } i_1 = 1 \\
                  p_{i_1,j_1}\prod_{k=2}^Np_{i_k,j_{N-k+2}},   \quad \text{if } i_1 = 0
                          \end{aligned} \right.
\end{equation}
where $\vec{i}=[i_1,i_2,\cdots,i_N]$, $\vec{j}=[j_1,j_2,\cdots,j_N]$, they are two ordered channel state vectors with entries equal to 0 or 1.

 Denote the probability vector of Markov chain formed by $\vec{S(t)}$ as $\vec{V}(t) = [v_1(t), v_2(t),\cdots, v_{2^N}(t)]$. Then we have
\begin{equation}
\vec{V}(t) = \vec{V}(1) \cdot (\mathbf{Q})^{t-1}
\end{equation}
In the $t_{th}$ step, the myopic policy selects the channel of the first component in $\vec{S(t)}$, therefore we only have to calculate the probability of states whose first component is 1. There are $2^{N-1}$ such states and denote them as $l_1, l_2, \cdots, l_{2^{N-1}} $, we have  \begin{equation}\label{sum}
\mathbb{E}[Y^{\pi_1}(\mathbf{P},\Omega(1),t)]
= \mathop{\Sigma}_{i=1}^{2^{N-1}}v_{l_i}(t)
\end{equation}

We can diagonalize $\mathbf{Q}$ as below:
\begin{equation}
\mathbf{Q} = \mathbf{H}^{-1}diag(\lambda_1,\cdots,\lambda_{2^N})\mathbf{H}
\end{equation}
From Perron Frobenius theorem \cite{Meyer:2000}, we know that
$|\lambda_i| \leq 1, \forall i$ and at least one eigenvalue is 1.
Without loss of generality, we assume $1=\lambda_1 = \cdots =
\lambda_a > |\lambda_{a+1}| \geq \cdots \geq |\lambda_{2^N}|,$ with
(\ref{sum}), we can rewrite $E[Y^{\pi_1}(\mathbf{P},\Omega(1),t)]$
as
\begin{equation}
\mathbb{E}[Y^{\pi_1}(\mathbf{P},\Omega(1),t)]
= \mathop{\Sigma}_{j=1}^{2^{N-1}}h_j\lambda_j^t
= \mathop{\Sigma}_{j=1}^{a}h_j + \mathop{\Sigma}_{j=a+1}^{2^{N-1}}h_j\lambda_j^t
\end{equation}
where $h_j$ is the corresponding coefficient which is only related to $\mathbf{P}$.

The steady average throughput $U_i(\mathbf{P})(i=1,2)$ is $\mathop{\Sigma}_{i=1}^{a}h_j$. As for different policies, the transition matrix $\mathbf{Q}$ is different, thus $\lambda$ vector and coefficient $h_j$ have different expressions. We denote them as $U_1(\mathbf{P})$ and $U_2(\mathbf{P})$ respectively.

Based on the same reason, for different policies, we have different
expressions for
$\mathop{\Sigma}_{j=a+1}^{2^{N-1}}|h_j|\frac{|\lambda_j|}{1-|\lambda_j|}$,
each denoted as $C_1(\mathbf{P})$ and $C_2(\mathbf{P})$
respectively.

Besides, we have:
\begin{equation}
\begin{split}
&|\mathop{\Sigma}_{t=1}^M{\Sigma}_{j=a+1}^{2^{N-1}}h_j\lambda_j^t | \leq \mathop{\Sigma}_{t=1}^M{\Sigma}_{j=a+1}^{2^{N-1}}|h_j||\lambda_j|^t\\
 &= \mathop{\Sigma}_{j=a+1}^{2^{N-1}}|h_j|{\Sigma}_{t=1}^M |\lambda_j|^t \leq \mathop{\Sigma}_{j=a+1}^{2^{N-1}}|h_j|{\Sigma}_{t=1}^\infty |\lambda_j|^t \\&\leq    \mathop{\Sigma}_{j=a+1}^{2^{N-1}}|h_j|\frac{|\lambda_j|}{1-|\lambda_j|}
\end{split}
\end{equation}

So $|\mathbb{E}^{\pi_i}[\Sigma_{t=1}^{M}Y^{\pi_i}(\mathbf{P},\Omega(1),t)]-M\cdot U_i(\mathbf{P})|<C_i(\mathbf{P})$, $i=1,2$.
\end{proof}

\begin{IEEEproof}[Proof of Theorem \ref{theorem:regret}]

We first derive a bound on the regret for the case when
$p_{01}<p_{11}$. In this case, policy $\pi_1$ would be the optimal.
Based on Lemma \ref{lemma:2}, the difference of
$\mathbb{E}^{\pi_1}[\Sigma_{t=1}^{n}Y^{\pi_1}(\mathbf{P},\omega_1,\omega_2,t)]$
and $nU_1$ is no more than $C_1$, therefore we only need to prove:
\begin{equation}
\begin{split}
&R'(\mathbf{P},\omega_1,\omega_2,n) \triangleq
nU_1-\mathbb{E}^{\pi_1}[\Sigma_{t=1}^{n}Y_{\pi_1}(\mathbf{P},\omega_1,\omega_2,t)]\\\nonumber
&\leq Z_1G(n)\ln(n)+Z_2\ln(n)+Z_3G(n)+Z_4-C_1
\end{split}
\end{equation}
where $Z_1,Z_2,Z_3,Z_4$ are
constants only related to $\mathbf{P}$.

The regret comes from two parts: the regret when using policy $\pi_2$; the regret between $U_1$ and
$Y^{\pi_1}(\mathbf{P},\omega_1,\omega_2,t)$ when using policy $\pi_1$. From Lemma \ref{lemma:2}, we know that each time when we switch
from policy $\pi_1$ to policy $\pi_2$, at most we lose a
constant-level value from the second part. So if the times of
policy $\pi_2$ being used is bounded by $O(G(n)\ln{n})$, both parts
of the rewards can be bounded by $O(G(n)\ln{n})$.

For case of exposition, we discuss the slots $n$ such that $G||n$, where $G||n$ denotes that time $n$ is the end of successive $G(n)$
plays.

We define $q$ as the smallest index such that
\begin{equation}
K_q \geq  \max \{
\lceil\frac{C_1+C_2}{|U_1-U_2|} \rceil, C_2/U_2, C_1/U_1 \}
\end{equation}

Let
$c_{t,s} = \sqrt{(L\ln{t})/s}$
,
$w_1 = q(U_1-\frac{C_1}{K_q})$
and
\begin{equation}
w_2 = q\frac{U_2-C_2/K_q}{U_2+C_2/K_q}(U_2+\frac{C_2}{K_q}-1)\nonumber
\end{equation}

 Next we will show that it is possible to
define $\alpha(U_1, C_1, \mathbf{P})$ such that if policy $\pi_1$
is played $s > \alpha$ times,
\begin{equation} \label{eqn:p1}
\exp(-2(w_1-sc_{t,s})^2/(s-q))
\leq t^{-4}.
\end{equation}

In fact, we have
\begin{equation}
\sqrt{Ls}-w_1 \geq \sqrt{2(s-q)}\nonumber
 \end{equation}
 when $s > \max{\{q, \lceil w_1/(\sqrt{L}-\sqrt{2})\rceil^2}\}$

 Consider
 \begin{equation}
 f(t) = \sqrt{Ls\ln{t}}-w_1- \sqrt{2(s-q)\ln{t}}\quad \forall t \geq e \nonumber
 \end{equation}

   It is quite obvious that $f(t)$ is a increasing function. Since $f(e) \geq 0$, we have $f(t) \geq 0, \forall t \geq e$, i.e.
  \begin{equation}
   \sqrt{Ls\ln{t}}-w_1\geq \sqrt{2(s-q)\ln{t}}   \nonumber
  \end{equation}
which equals to
\begin{equation}
\exp(-2(w_1-sc_{t,s})^2/(s-q))\leq t^{-4}\nonumber
\end{equation}
 Thus at least we can set
 \begin{equation}
 \alpha(U_1, C_1, \mathbf{P}) = \max{\{q, \lceil w_1/(\sqrt{L}-\sqrt{2})\rceil^2}\} \nonumber
 \end{equation}

For the similar reason, we could define
\begin{equation}
\beta(U_2, C_2, \mathbf{P}) =  \max{\{q,\lceil w_2/(\sqrt{L}-\sqrt{2})\rceil ^2}\}   \nonumber
\end{equation}
 such that if
policy $\pi_2$ is played $s  > \beta$ times,
\begin{equation}\label{eqn:p2}
\exp(\frac{-2(w_2+sc_{t,s})^2}{s-q})\\
\leq t^{-4}
\end{equation}

Moreover, we will show that there exists

\noindent $\gamma = \lceil
\max\{{5\alpha+1,e^{4\alpha/L}+\alpha,5\beta+1,e^{4\beta/L}+\beta}\}\rceil$
such that when $G(n)>K_{\gamma}$, policy $\pi_1$ is played at least
$\alpha$ times and policy $\pi_2$ is played at least $\beta$ times.

In fact, if policy $\pi_1$ has been played less than $\alpha$ times, then policy $\pi_2$ has been played at least $(4\alpha+2)$ times. Consider the last time selecting policy $\pi_2$, there must be
\begin{equation}
\frac{\hat{X}_{1,i_1}}{i_1}+c_{t,i_1}\leq\frac{\hat{X}_{2,i_2}}{i_2}+c_{t,i_2}
\end{equation}
Noting that $\frac{\hat{X}_{1,i_1}}{i_1} \geq 0$, $\frac{\hat{X}_{2,i_2}}{i_2} \leq 1$, $i_1 \leq \alpha -1$,$i_2 \geq 4\alpha + 1$,  we have
\begin{equation}
0 + \sqrt{\frac{L\ln{t}}{\alpha-1}} \leq 1+ \sqrt{\frac{L\ln{t}}{4\alpha+1}}    \nonumber
\end{equation}

Consider
\begin{equation}
g(t) =  1+ \sqrt{\frac{L\ln{t}}{4\alpha+1}}- \sqrt{\frac{L\ln{t}}{\alpha-1}}     \nonumber
\end{equation}
 Since $g(t)$ is a decreasing function and $t \geq \gamma-\alpha  \geq e^{4\alpha/L}$, we have
 \begin{equation}
 g(t) \leq g(e^{4\alpha/L}) = 1+ \sqrt{\frac{4\alpha}{4\alpha+1}}- \sqrt{\frac{4\alpha}{\alpha-1}} < 0  \nonumber
 \end{equation}
which contradicts the conclusion above. So policy $\pi_1$ has been played at least $\alpha$ times. For the similar reason, policy $\pi_2$ is played at least $\beta$ times.

Denote $T(n)$ as the number of times we select policy $\pi_2$ up to
time $n$. Then, for any positive integer $l$, we have
\begin{equation}
\begin{split} T(n) & =
1+\sum_{t=K_1+K_2,G||t}^{n}\mathbb{I}\{\frac{\hat{X}_1(t)}{i_1(t)}+c_{t,i_1(t)} <
\frac{\hat{X}_2(t)}{i_2(t)}+\\
&c_{t,i_2(t)}\} \\
& \leq l+\gamma+\\
&\sum_{t=K_1+\cdots+K_\gamma,G||t}^{n}\sum_{s_1=\alpha}^{\alpha(t),t=K_1+\cdots+K_{\alpha(t)}}\sum_{s_2=max(\beta,l)}^{\beta(t),t=K_1+\cdots+K_{\beta(t)}}\\
&\mathbb{I}\{\frac{\hat{X}_{1,s_1}}{s_1}+c_{t,s_1}\leq\frac{\hat{X}_{2,s_2}}{s_2}+c_{t,s_2}\}
\end{split}
\end{equation}
where $\mathbb{I}\{x\}$ is the index function defined to be 1 when
the predicate $x$ is true, and 0 when it is false predicate;
$i_j(t)$ is the number of times we select policy $\pi_j$ when up to
time $t, \forall j=1,2$; $\hat{X}_j(t)$ is the sum of every sample
mean for $K_i$ plays up to time $t$; $\hat{X}_{i,s_i}$ is the sum of
every sample mean for $K_{s_i}$ times using policy $\pi_i$.

The condition
$\{\frac{\hat{X}_{1,s_1}}{s_1}+c_{t,s_1}\leq\frac{\hat{X}_{2,s_2}}{s_2}+c_{t,s_2}\}$
implies that at least one of the following must hold:
\begin{equation} \label{eqn:inequ1}
\frac{\hat{X}_{1,s_1}}{s_1}\leq
U_1-\frac{C_1}{K_q}-c_{t,s_1}
\end{equation}
\begin{equation} \label{eqn:inequ2}
\frac{\hat{X}_{2,s_2}}{s_2}\geq
U_2+\frac{C_2}{K_q}+\frac{U_2+C_2/K_q}{U_2-C_2/K_q}c_{t,s_2}
\end{equation}
\begin{equation} \label{eqn:inequ3}
U_1-\frac{C_1}{K_q}<U_2+\frac{C_2}{K_q}+(1+\frac{U_2+C_2/K_q}{U_2-C_2/K_q})c_{t,s_2}
\end{equation}

Note that
\begin{equation}
\hat{X}_{1,s_1}=\hat{A}_{1,1}+\hat{A}_{1,2}+\cdots+\hat{A}_{1,s_1}
\end{equation}
where $\hat{A}_{1,i}$ is sample average reward for the $i_{th}$ time
selecting policy $\pi_1$.

Due to the definition of $\alpha$ and $K_q$, we have
\begin{equation}
U_1-\frac{C_1}{K_q} \leq \mathbb{E}[\hat{A}_{1,i}]\leq U_1+\frac{C_1}{K_q} \quad \forall i \geq q
\end{equation}

 Then applying Lemma \ref{lemma:chernoff}, and the
results in (\ref{eqn:p1}) and (\ref{eqn:p2}),
\begin{equation}
\begin{split}
&\mathbb{P}(\frac{\hat{X}_{1,s_1}}{s_1}\leq
U_1-\frac{C_1}{K_q}-c_{t,s_1}) \\
&=\mathbb{P}(\frac{\hat{A}_{1,1}+\hat{A}_{1,2}+\cdots+\hat{A}_{1,s_1}}{s_1}\leq
U_1-\frac{C_1}{K_q}-c_{t,s_1}) \\
&\leq \mathbb{P}(\frac{0+\cdots+0+\hat{A}_{1,q+1}+\hat{A}_{1,2}+\cdots+\hat{A}_{1,s_1}}{s_1}\leq
U_1\\
&-\frac{C_1}{K_q}-c_{t,s_1})\\
&\leq \exp(-2(w_1-s_1c_{t,s_1})^2/(s_1-q))\leq t^{-4}
\end{split}
\end{equation}

Similarly,

\begin{equation}
\begin{split}
&\mathbb{P}(\frac{\hat{X}_{2,s_2}}{s_2}\geq
U_2+\frac{C_2}{K_q}+\frac{U_2+C_2/K_q}{U_2-C_2/K_q}c_{t,s_2})\\
&=\mathbb{P}(\frac{\hat{A}_{2,1}+\hat{A}_{2,2}+\cdots+\hat{A}_{2,s_2}}{s_2}\geq
U_2+\frac{C_2}{K_q}\\
&+\frac{U_2+C_2/K_q}{U_2-C_2/K_q}c_{t,s_2})\\
&\leq \mathbb{P}(\frac{1+\cdots+1+\hat{A}_{2,q+1}+\cdots+\hat{A}_{2,s_2}}{s_2}\geq
U_2\\
&+\frac{C_2}{K_q}+\frac{U_2+C_2/K_q}{U_2-C_2/K_q}c_{t,s_2})\\
&\leq \exp(\frac{-2(w_2+s_2c_{t,s_2})^2}{s_2-q})\\
&\leq t^{-4}
\end{split}
\end{equation}


Denote $\lambda(n)$ as
\begin{equation}
\lambda(n)=\lceil(L(1+\frac{U_2+C_2/K_q}{U_2-C_2/K_q})^2\ln{n})/(U_1-U_2-\frac{C_1+C_2}{K_q})^2\rceil
\end{equation}

For $l \geq \lambda(n)$, (\ref{eqn:inequ3}) is false. So we get:
\begin{equation}
\begin{split}
&\mathbb{E}(T(n)) \leq
\lambda(n) +\gamma+
\Sigma_{t=1}^{\infty}\Sigma_{s_1=1}^{t}\Sigma_{s_2=1}^{t}2t^{-4}\\
&\leq
\lambda(n)+\gamma+\frac{\pi^2}{3}.
\end{split}
\end{equation}

Therefore, we have:
\begin{equation}
\begin{split}
& R'(\mathbf{P},\Omega(1),n) \leq G(n)+ \\
& \quad\quad  ((U_1-U_2)G(n)+C_1+C_2+C_1)
(\lambda(n)+\gamma+\frac{\pi^2}{3})
\end{split}
\end{equation}

This concludes the bound in case $p_{11}>p_{01}$. The derivation of
the bound is similar for the case when $p_{11}\leq p_{01}$ with the
key difference of $ \gamma'$ instead of $\gamma$, and the $C_1, U_1$
terms being replaced by $C_2, U_2$ and vice versa. Then we have that
the regret in either case has the following bound:

\begin{equation}
\begin{split}
&R(\mathbf{P},\Omega(1),n) \leq G(n)+
(|U_2-U_1|G(n)+C_1+C_2+\\
&\max\{C_1,C_2\})
(\frac{L(1+\max\{\frac{U_1+C_1/K_q}{U_1-C_1/K_q},\frac{U_2+C_2/K_q}{U_2-C_2/K_q}\})^2\ln{n}}{(|U_2-U_1|-\frac{C_1+C_2}{K_q})^2}
\\
& + 1 +{\max\{\gamma,\gamma'}\}+\frac{\pi^2}{3})+\max\{C_1,C_2\}
\end{split}
\end{equation}

This inequality can be readily translated to the simplified form of
the bound given in the statement of Theorem 1, where:
\begin{equation}
\begin{split}
& Z_1 = |U_2-U_1|
\frac{L(1+\max\{\frac{U_1+C_1/K_q}{U_1-C_1/K_q},\frac{U_2+C_2/K_q}{U_2-C_2/K_q}\})^2}{(|U_2-U_1|-\frac{C_1+C_2}{K_q})^2}\\
& Z_2 = (C_1+C_2+\max\{C_1,C_2\})
L(1+\\
&\max\{\frac{U_1+C_1/K_q}{U_1-C_1/K_q},\frac{U_2+C_2/K_q}{U_2-C_2/K_q}\})^2/(|U_2-U_1|-\frac{C_1+C_2}{K_q})^2\\
& Z_3 = |U_2-U_1| (1 +{\max\{\gamma,\gamma'}\}+\frac{\pi^2}{3}) +
1\\ \nonumber
& Z_4 = (C_1+C_2+\max\{C_1,C_2\}) (1
+{\max\{\gamma,\gamma'}\}+\frac{\pi^2}{3})\\
&+\max\{C_1,C_2\}
\end{split}
\end{equation}
\end{IEEEproof}

\textbf{Remark 1:} Theorem 1 has been stated for the cases $N=2,3$,
which are the cases when the Myopic policy has been proved to be
optimal for the known-parameter case for all values of $P$. However,
our proof shows something stronger: Algorithm \ref{alg:sensing}
yields the claimed near-logarithmic regret with respect to the
myopic policy for any $N$. The Myopic policy is known to be always
optimal for $N=2,3$, and for any $N$ so long as the Markov chain is
positively correlated. For negatively correlated channels, it is an
open question whether it is optimal for an infinite horizon case
(extensive numerical examples suggest an affirmative answer to this
conjecture). If this conjecture is true, the policy we have
presented would also offer near-logarithmic regret for any $N$.

\textbf{Remark 2:} Theorem 1 also holds if the initial belief vector
is unknown. This is because once every channel is sensed once,
initial belief is forgotten by the myopic policy, which must happen
within finite time on average.

\section{Examples And Simulation Results}

We consider a system that consists of $N = 2$ independent channels.
Each channel evolves as a two-state Markov chain with transition
probability matrix $P$. The parameter $L$ is set to be 3. We
consider several situations with different sequence
$\{K_n\}_{n=1}^\infty$ and different correlations.

First we show the simulation results when the channel is positively
correlated. The transition probabilities are as follows: $p_{01} =
0.3$, $p_{11} = 0.7$, $p_{10} = 0.3$, $p_{00} = 0.7$. The sequences
are set to be $K^1_n=\lceil 100+\ln(n+2)\rceil$, $K^2_n=\lceil
100+\ln(\ln(n+2))\rceil$ and $K^3_n=\lceil
100+\ln(\ln(\ln(n+2)))\rceil$.

\begin{figure}[ht]
 \centering
 \includegraphics[width=0.5\textwidth]{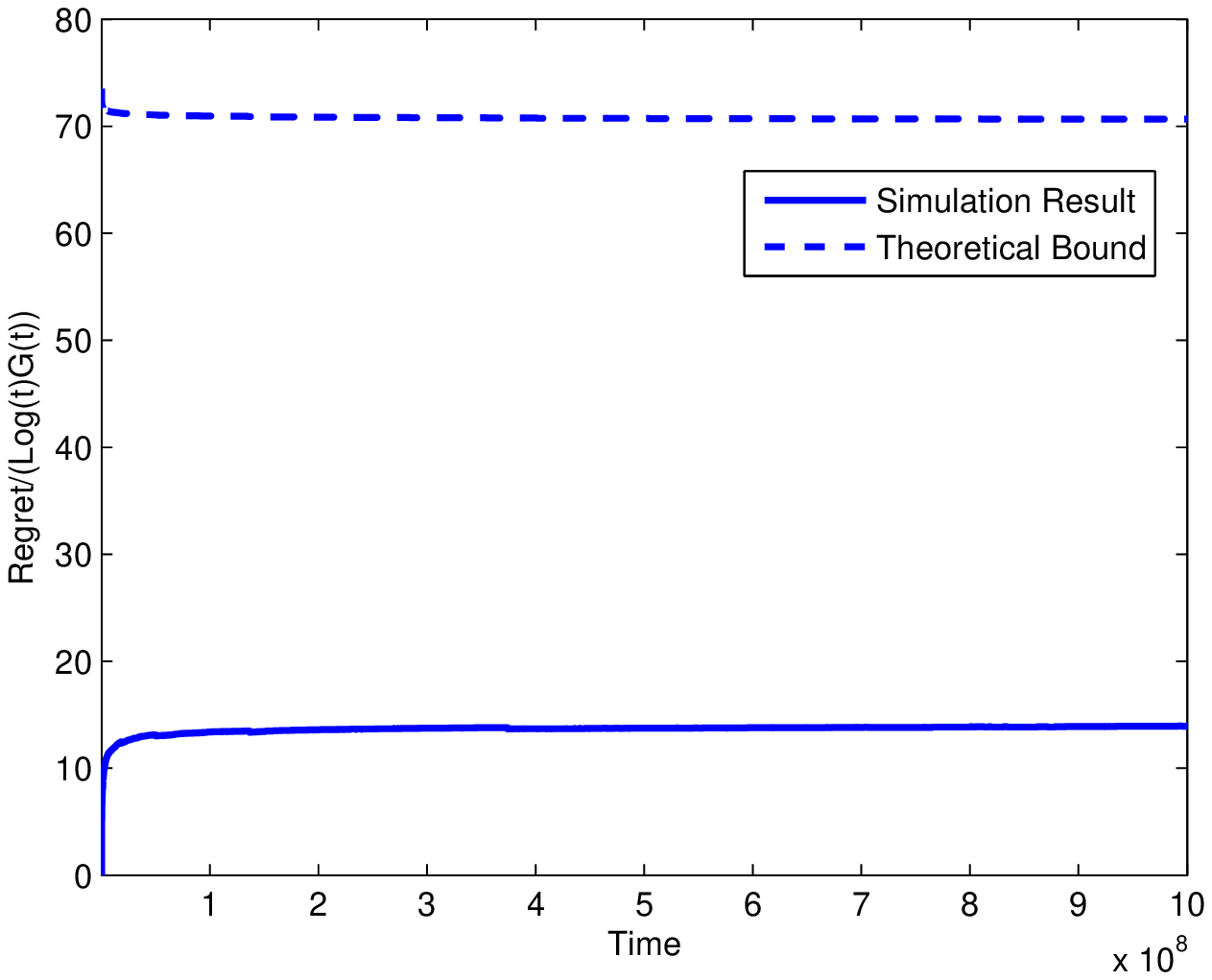}
\caption{Simulation results when $K_n = \lceil 100+\ln(n+2)\rceil$} \label{fig:2}
\end{figure}

\begin{figure}[ht]
 \centering
 \includegraphics[width=0.5\textwidth]{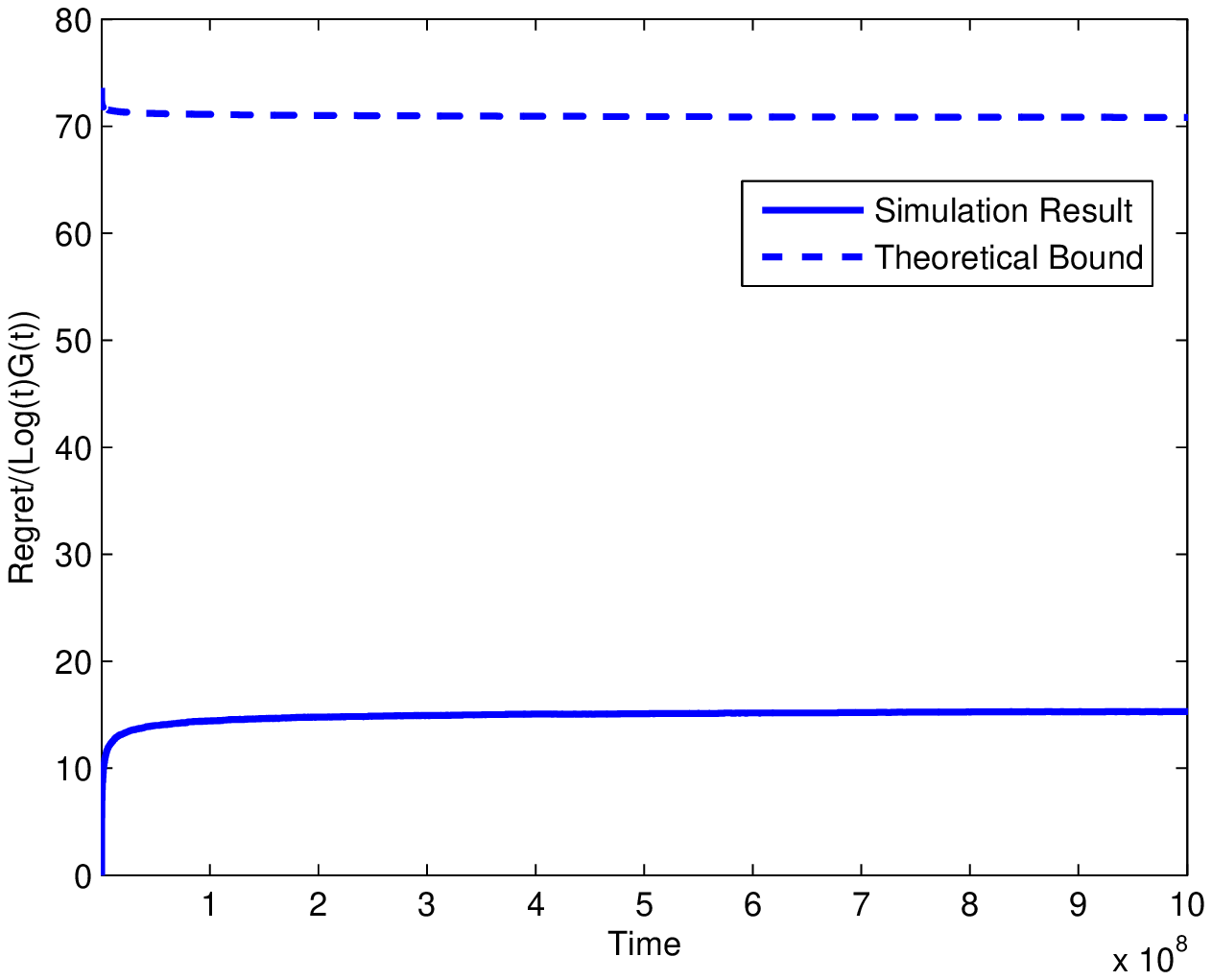}
\caption{Simulation results when $K_n = \lceil 100+\ln(\ln(n+2))\rceil$} \label{fig:1}
\end{figure}

\begin{figure}[ht]
 \centering
 \includegraphics[width=0.5\textwidth]{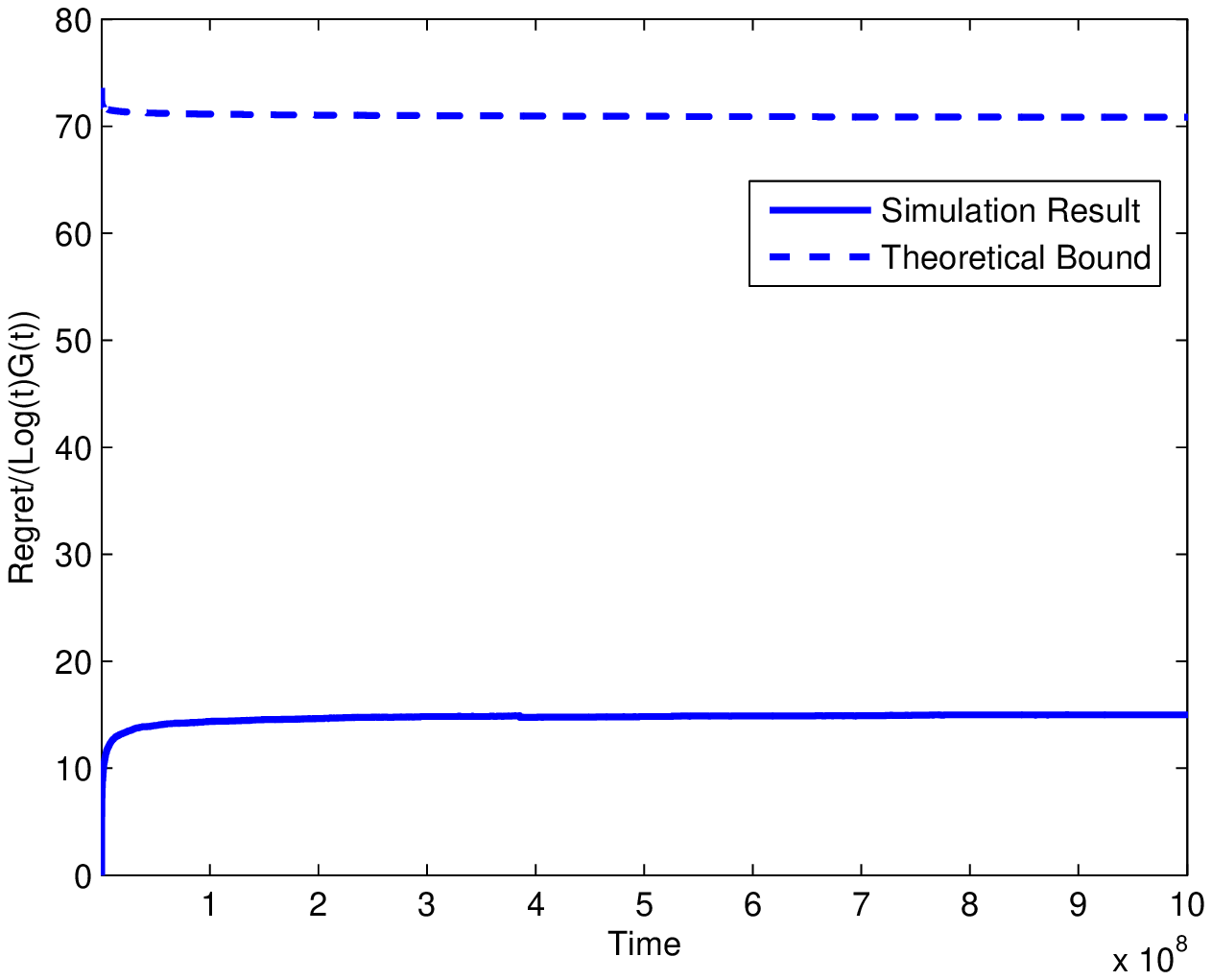}
\caption{Simulation results when $K_n = \lceil 100+\ln(\ln(\ln(n+2))\rceil$)} \label{fig:3}
\end{figure}

Figure \ref{fig:1} to Figure \ref{fig:3} show the simulation result (normalized with respect to $G(t)\log t$. It is quite clear that the regrets all converge
to a limit that is below our bound. In our simulations, $\lceil
\frac{C_1+C_2}{|U_1-U_2|} \rceil$ is 6 and $K_1$ is already greater
than it. In this way, the regret can converge more quickly.
Practically, it may happen that $K_1 < \lceil
\frac{C_1+C_2}{|U_1-U_2|} \rceil$. Then, we have to wait for some
time so $K_n$ can be sufficiently great. Since $K_n$ goes to the
infinity, it will exceed $\lceil \frac{C_1+C_2}{|U_1-U_2|} \rceil$
at certain time. The speed of convergence depends on how fast $K_n$
grows. If $K_n$ grows too slowly, it may take longer time to
converge; however, if it grows too fast, though the regret converges
quickly, the upper bound of regret also increases. So there is
trade-off here between convergence speed and the upper bound of
regret. Generally $K_n$ should be a sub-linear sequence.

Next, we show the simulation results when the channel is negatively
correlated. The transition probabilities are as follows: $p_{01} =
0.7$, $p_{11} = 0.3$, $p_{10} = 0.7$, $p_{00} = 0.3$. We use again
the sequence $K^1_n$, $K^2_n$ and $K^3_n$.

 \begin{figure}[ht]
 \centering
 \includegraphics[width=0.5\textwidth]{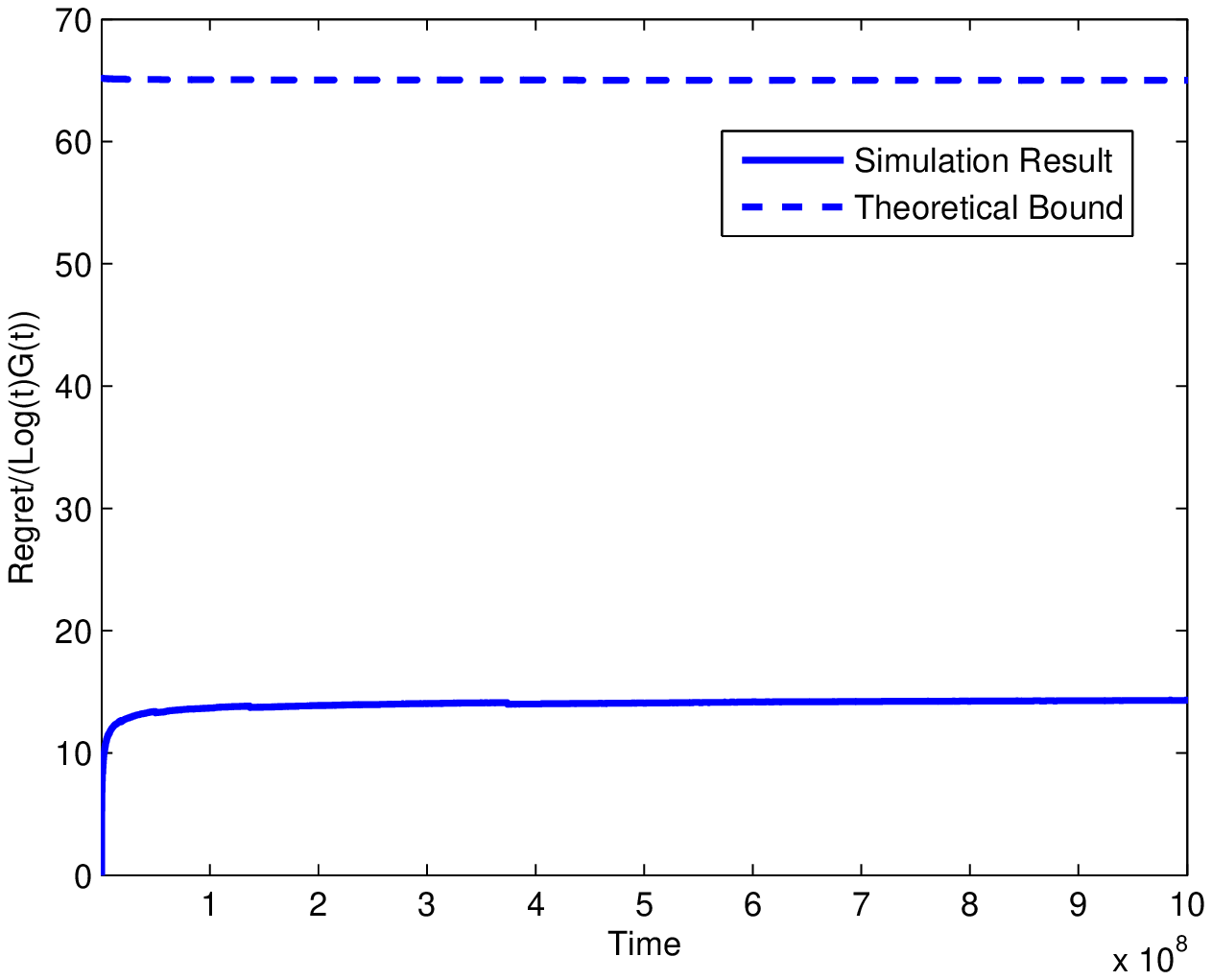}
\caption{Simulation results when $K_n = \lceil 100+\ln(n+2)\rceil$} \label{fig:4}
\end{figure}

\begin{figure}[ht]
 \centering
 \includegraphics[width=0.5\textwidth]{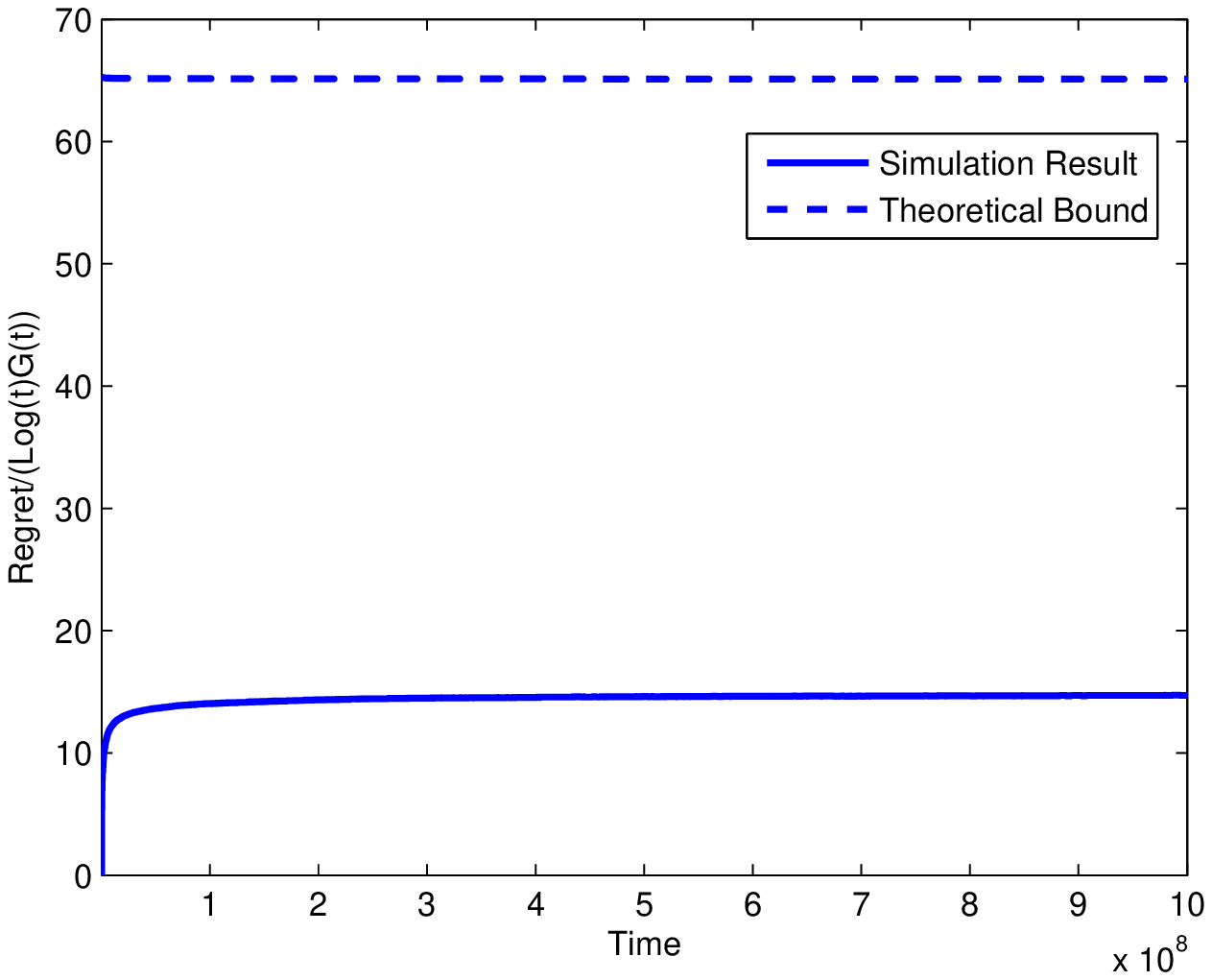}
\caption{Simulation results when $K_n = \lceil 100+\ln(\ln(n+2))\rceil$} \label{fig:5}
\end{figure}

\begin{figure}[ht]
 \centering
 \includegraphics[width=0.5\textwidth]{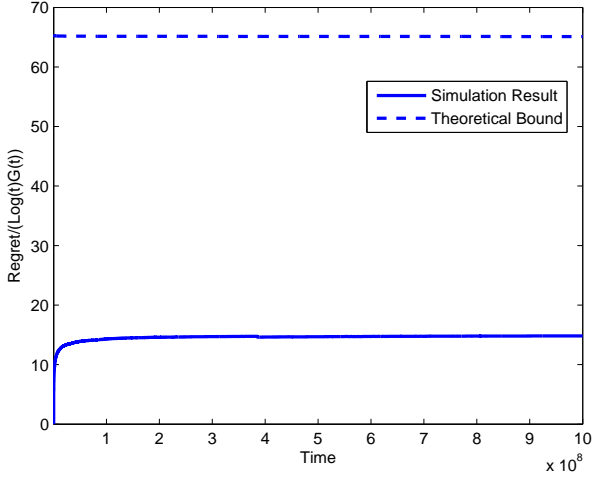}
\caption{Simulation results when $K_n = \lceil 100+\ln(\ln(\ln(n+2))\rceil$)} \label{fig:6}
\end{figure}

Figure \ref{fig:4} to Figure \ref{fig:6} show the simulation of the
result (normalized with respect to the product of $G(t)$ and the
logarithm of time). The regrets also converge to a limitation and
are bounded. The basic conclusion also stands here.

\section{Conclusion}
In this study, we have considered the non-Bayesian RMAB problem, in
which the parameters of Markov chain are unknown \emph{a priori}. We
have developed a novel approach to solves special cases of this
problem that applies whenever the corresponding known-parameter
Bayesian problem has the structure that the optimal solution is one
of the prescribed finite set of polices depending on the known
parameters. For such problems, we propose to learn the optimal
policy by using a meta-policy which treats each policy from that
finite set as an arm. We have demonstrated our approach by
developing an original policy for opportunistic spectrum access over
unknown dynamic channels. We have proved that our policy achieves
near-logarithmic regret in this case, which is the first-such strong
regret result in the literature for a non-Bayesian RMAB problem.

While we have demonstrated this meta-policy approach for a
particular RMAB with two states and identical arms, an open question
is to identify other RMAB problems that fill into the finite option
structure, and derive similar results for them. Note that even for
the problem where the optimal solution does not have a finite option
structure, but there exists a near-optimal policy that has this
structure, our approach could be used to prove sub-linear regret
with respect to the near-optimal policy.

%
%

\ifCLASSOPTIONcaptionsoff
  \newpage
\fi

\appendices
\section{Calculation of $C_i({\mathbf{P}})$ and $U_i({\mathbf{P}})$ for Lemma \ref{lemma:2} when $N = 2$}\label{appendix:01}

When $N=2$, we explicitly calculate $C_1({\mathbf{P}})$,
$C_2({\mathbf{P}})$, $U_1({\mathbf{P}})$, $U_2({\mathbf{P}})$ as
follows:
\begin{equation}\label{eq:w}
\begin{split}
U_1({\mathbf{P}}) = &
\frac{p_{01}^2}{(1-p_{11}+p_{01})^2}+\frac{(1-p_{01}+p_{11})p_{01}(1-p_{11})}{(1-p_{11}+p_{01})^2}
\end{split}
\end{equation}

\begin{equation}
\begin{split}
U_2({\mathbf{P}}) = &
\frac{p_{01}^2}{(1-p_{11}+p_{01})^2}+\frac{p_{01}(1-p_{11})}{1-p_{11}+p_{01}}
\end{split}
\end{equation}

\begin{equation}
\begin{split}
C_1({\mathbf{P}}) & = \frac{2\max\{p_{01},1-p_{11}\}|p_{11}-p_{01}|(1-p_{11})}{(1-p_{11}+p_{01})^3} \\
&
+\frac{\max\{p_{01}^2,(1-p_{11})^2\}|p_{11}-p_{01}|}{(1-(p_{11}-p_{01})^2)(1-p_{11}+p_{01})^2}
\end{split}
\end{equation}

\begin{equation}
\begin{split}
C_2({\mathbf{P}}) & = \frac{2\max\{p_{01},1-p_{11}\}|p_{11}-p_{01}|p_{01}}{(1-p_{11}+p_{01})^3} \\
&
+\frac{\max\{p_{01}^2,(1-p_{11})^2\}|p_{11}-p_{01}|}{(1-(p_{11}-p_{01})^2)(1-p_{11}+p_{01})^2}
\end{split}
\end{equation}

Also, note that a similar but more tedious calculation can be done
for $N = 3$.


\end{document}